\documentclass{amsart}

\usepackage{amssymb,amsmath,latexsym,amsfonts,enumerate,calc,verbatim}

\newcommand\C{{\mathbb C}}

\newcommand\R{{\mathbb R}}

\newtheorem{theorem}{Theorem}

\newtheorem{lemma}{Lemma}
\newtheorem{cor}{Corollary}
\newtheorem{proposition}{Proposition}

\begin{document}

\title[Maps conjugating holomorphic maps]{Maps conjugating holomorphic
  maps in $\C^n$}
\author[G.~Buzzard, S.~Merenkov]{Gregery~T.~Buzzard, and Sergei~Merenkov}

\address{Department of Mathematics\\ Purdue University\\ West
  Lafayette, IN 47907-2067}

\email{buzzard@math.purdue.edu}

\address{Department of Mathematics\\ The University of Michigan\\
Ann Arbor, MI 48109-1109}

\email{merenkov@umich.edu}

\abstract {If $\psi$ is a bijection from $\C^n$ onto a complex
manifold $\mathcal M$, which conjugates every holomorphic map in
$\C^n$ to an endomorphism in $\mathcal M$, then we prove that
$\psi$ is necessarily biholomorphic or antibiholomorphic. This
extends a result of A.~Hinkkanen to higher dimensions. As a
corollary, we prove that if there is an epimorphism from the
semigroup of all holomorphic endomorphisms of $\C^n$ to the
semigroup of holomorphic endomorphisms of a complex manifold
$\mathcal M$ consisting of more than one point, or an epimorphism
in the opposite direction for a doubly-transitive $\mathcal M$,
then it is given by conjugation by some biholomorphic or
antibiholomorphic map. We show also that there are two unbounded
domains in $\C^n$ with isomorphic endomorphism semigroups but
which are neither biholomorphically nor antibiholomorphically
equivalent.}
\endabstract

\maketitle


\section{Introduction}

The question of determining a mathematical structure of an object
from its semigroup of endomorphisms, i.e. the set of all maps from
this object into itself with composition as a semigroup operation,
goes back to at least C.~J.~Everett and S.~M.~Ulam~\cite{cE48},
\cite{sU60}. In the most general form this question can be
formulated as follows. Suppose we have two sets $A$ and $B$ with a
given structure, whose semigroups of endomorphisms compatible with
this structure are isomorphic. Does there exist a bijective map
between $A$ and $B$, which preserves the structure?

K.~D.~Magill, L.~M.~Glusk{\=\i}n, B.~M.~Schein,
L.~B.~{\v{S}}neperman, and I.~S.~Yoroker studied the question of
determining a topological space from its semigroup of continuous
endomorphisms. See a survey in~\cite{kM75}, \cite{lG77}.

To the best of the authors' knowledge, L.~Rubel was the first who
raised the question of determining a complex space from its
semigroup of holomorphic endomorphisms. In 1993,
A.~Eremenko~\cite{aE93} proved that every two Riemann surfaces
that admit non-constant bounded holomorphic functions, and whose
semigroups of holomorphic endomorphisms are isomorphic, are
necessarily conformally or anticonformally equivalent. This result
was extended by S.~Merenkov in~\cite{sM02} to bounded domains in
$\C^n$.

On the other hand, A.~Hinkkanen~\cite{aH92} proved in 1992 that
there exist unbounded domains in $\C$ whose semigroups of
holomorphic endomorphisms are isomorphic, but the domains are not
even homeomorphic. In the same paper A.~Hinkkanen studied another
question raised by L.~Rubel. Namely, he proved that if $\psi$ is a
one-to-one function of the plane onto itself (but with no
assumption of continuity), such that $\psi\circ f\circ\psi^{-1}$
is entire whenever $f$ is entire, then $\psi$ has the form
$\psi(z)=az+b$, or $\psi(z)=a\overline{z}+b$, where $a$ and $b$
are complex numbers with $a\neq 0$;  i.e., $\psi$ is a conformal
or anticonformal automorphism.

In higher dimensions, any analog of A.~Hinkkanen's theorem must
take into account the fact that the automorphism group of $\C^n$
is quite large, since in $\C^2$, for example, there are
biholomorphic maps of the form $\psi(z_1, z_2)=(z_1, z_2+g(z_1))$,
where $g$ is an arbitrary entire function. However, one may hope
that arbitrary $\psi$ conjugating holomorphic maps to holomorphic
maps is still a biholomorphic or antibiholomorphic automorphism.
The main theorem, Theorem~\ref{T:1}, of the present paper asserts
that this is indeed the case.

Note that the set of all holomorphic endomorphisms of a complex
manifold $\mathcal M$ forms a semigroup (with unit) under
composition. We denote this semigroup by $E({\mathcal M})$. If
${\mathcal M}=\C^n$, we denote the semigroup by $E$.

\begin{theorem}\label{T:1}
If $\psi$ is a bijection of $\C^n,\ n\geq 2$ onto a
complex manifold
$\mathcal M$, such that $\psi\circ f\circ\psi^{-1}\in E({\mathcal M})$
for every map $f\in E$, then $\psi$ is biholomorphic or
antibiholomorphic.
\end{theorem}

As in the one-dimensional case~\cite{aH92}, it is not sufficient
to assume that $\psi\circ f\circ\psi^{-1}\in E({\mathcal M})$ for every
polynomial map $f$ in order to conclude that $\psi$ is a homeomorphism.
The reason is that there are non-continuous field automorphisms of
$\C$~\cite{sL72}. If $\xi$ is such an automorphism, then we can
take $\psi(z_1,\dots, z_n)=(\xi(z_1),\dots, \xi(z_n))$. The
conjugation by $\psi$ is an automorphism of semigroups of
polynomial maps in $\C^n$, but $\psi$ is not continuous.

We say that a complex manifold $\mathcal N$ is
{\emph{doubly-transitive}} if $E({\mathcal N})$ is
doubly-transitive, i.e. if for every pair $z_1, z_2$ of distinct
points in $\mathcal N$ and every other pair of points $w_1, w_2$
in $\mathcal N$, there exists $f\in E({\mathcal N})$, such that
$f(z_m)=w_m,\ m=1, 2$. We say that $\mathcal N$ is {\emph{weakly
doubly-transitive}} if in the previous definition we replace the
assumption that $w_2$ arbitrary by requiring that it has to be
sufficiently close to $w_1$. Clearly, every doubly-transitive
complex manifold is weakly doubly-transitive, and $\C^n$ is
doubly-transitive. As a corollary to Theorem~\ref{T:1}, we prove
the following
\begin{theorem}\label{T:2}
If there exists an epimorphism
of semigroups  $\phi:\ E\to E({\mathcal M})$, where $\mathcal M$ is
a complex manifold consisting of more than one point, then
\begin{equation}\label{E:M}
\phi(f)=\psi\circ f\circ\psi^{-1},\ \ \forall f\in E,
\end{equation}
for some biholomorphic or antibiholomorphic map $\psi:\ \C^n\to
{\mathcal M}$.

If there exists an epimorphism of semigroups $\varphi:\
E({\mathcal M})\to E$, where $\mathcal M$ is a weakly
doubly-transitive complex manifold, then
\begin{equation}\label{E:M2}
\varphi(f)=\eta\circ f\circ\eta^{-1},\ \ \forall f\in E({\mathcal
M}) ,
\end{equation}
for some biholomorphic or antibiholomorphic map $\eta:\ {\mathcal
M} \to \C^n$.
\end{theorem}

We note that the converse to this theorem is
trivial. If $\psi$ is a biholomorphic  or antibiholomorphic map from
$\C^n$ to $\mathcal M$, then the map $f\mapsto \psi\circ
f\circ\psi^{-1}$ is an isomorphism between the semigroups.  Similarly,
we get an isomorphism of semigroups if there exists an
(anti)biholomorphic map $\eta:\  {\mathcal M}\to \C^n$.  In
particular, we obtain the following corollary, which follows immediately
from the previous remarks plus the fact that an antibiholomorphic
equivalence from $\C^n$ to $\mathcal M$ implies a biholomorphic
equivalence simply by composing with the involution $z \mapsto
\overline{z}$.

\begin{cor} \label{C:1}
Given a complex manifold $\mathcal M$, the endomorphism semigroup of
$\mathcal M$ is isomorphic to the endomorphism semigroup of $\C^n$ if
and only if $\mathcal M$ is biholomorphic to $\C^n$.
\end{cor}

The first part of Theorem~\ref{T:2} is in some sense quite
surprising because, among the complex manifolds of dimension $n$,
$\C^n$ has a large and complicated semigroup of endomorphisms
(compare the simple semigroups in Theorem~\ref{T:3} below).  Yet
the equivalence given above requires only the existence of an
epimorphism from the ``large'' semigroup $E$ onto $E({\mathcal
M})$.

Also, applying methods used by D.~Varolin~\cite{dV00}, any Stein
manifold $\mathcal M$ with the (volume) density property is
doubly-transitive, and hence can be used in the second part of
Theorem~\ref{T:2}. Indeed, the fact that the manifold is Stein
implies that any single point is a holomorphically convex set.
Then for distinct points $p_1, p_2, q_1, q_2$ in $\mathcal M$,
Theorem~0.2 of~\cite{dV00} with $K = \{p_2\}$ implies that there
is an automorphism, $f_1$, of $\mathcal M$ so that $f_1(p_1) =
q_1$ and $f_1(p_2) = p_2$.  Likewise, there is an automorphism,
$f_2$, of $\mathcal M$ so that $f_2(p_2) = q_2$ and $f_2(q_1) =
q_1$. Thus for the map $f = f_2 \circ f_1$ we have $f(p_j) = q_j$.
If $p_1, p_2$ are distinct and $q_1, q_2$ are arbitrary (and
possibly one or both of them is the same as $p_1$ or $p_2$), then
we can first choose $z_1 \neq z_2$, distinct from the previous
four points, map $p_j$ to $z_j$, and then $z_j$ to $q_j$ (using a
constant map if $q_1 = q_2$).  Hence $\mathcal M$ is
doubly-transitive.

We mention also a recent paper by
S.~G.~Krantz~\cite{sK02}, where he studies the question of
determination of a domain in complex space by its automorphism
group. Of course a domain possesses more endomorphisms than
automorphisms. Therefore the ability to determine a domain
from its automorphism group implies the ability to determine a domain
from its endomorphism semigroup. Our Theorem~\ref{T:2} differs from
Krantz's result in that, first of all, we assume the
existence of an epimorphism between semigroups, rather than an
isomorphism. Secondly, the information we assume has a purely
algebraic character, i.e. the existence of an algebraic
epimorphism, and not a topological one; i.e., we make no a priori
assumptions about continuity.  To our knowledge it is an
open question if the existence of a purely algebraic isomorphism between the
automorphism group of $\C^n$ and the automorphism group of ${\mathcal
  M}$ implies the biholomorphic equivalence of these manifolds.
However, one result along these lines is contained
in the work of P.~Ahern and W.~Rudin~\cite{pA95}. They showed
that ${\text{Aut}}(\C^n)$ is sensitive to the dimension, i.e. if $1\leq m< n$,
then the groups ${\text{Aut}}(\C^m)$ and ${\text{Aut}}(\C^n)$
are not algebraically isomorphic.

To complete the analogy with Hinkkanen's results, we show the
existence of two unbounded domains in $\C^n$ with isomorphic
endomorphism semigroups but which are not (anti)bi\-holo\-morphic\-al\-ly
equivalent.  This should be compared with Merenkov's result
\cite{sM02}, in which it is shown that for two bounded domains in $\C^n$,
an isomorphism between the endomorphism semigroups implies the
(anti)bi\-holomorph\-ic equivalence between the two domains.

\begin{theorem} \label{T:3}
There exist unbounded domains $D_1$ and $D_2$ in $\C^n$ so that the
endomorphism semigroups $E(D_1)$ and $E(D_2)$ are isomorphic but such
that there is no biholomorphic or antibiholomorphic map from $D_1$
onto $D_2$.
\end{theorem}

The paper is organized as follows. In Section~\ref{S:FB} we prove
that the map $\psi$ in Theorem~\ref{T:1} is a homeomorphism, using the
notion of a Fatou-Bieberbach domain and pose a question about
Fatou-Bieberbach domains in Stein manifolds with the density
property. Section~\ref{S:LL} and Section~\ref{S:ME}
are devoted to the proof that $\psi$ is biholomorphic or
antibiholomorphic. In Section~\ref{S:E} we give a proof
of Theorem~\ref{T:2}, and in Section~\ref{S:I}, we prove
Theorem~\ref{T:3}.

{\bf{Acknowledgement.}} The authors thank A.~Eremenko, A.~Hinkkanen,
S.~Krantz, and J.P.~Rosay for helpful conversations.

\section{Fatou-Bieberbach domains and continuity of
$\psi$}\label{S:FB}

Below we assume that $n\geq 2$.

Let $FB$ denote the set of Fatou-Bieberbach domains, i.e. proper
domains in $\C^n$ that are biholomorphic to $\C^n$. A domain from
this set will be called an $FB$-{\emph{domain}}, and a
biholomorphic map from $\C^n$ onto an $FB$-domain will be called
an $FB$-{\emph map}.

We denote by $\Delta(r)$ the disk in $\C$ centered at 0 and of
radius $r$, and by $\Delta^k(r)$ the $k$-fold product of
$\Delta(r)$. In~\cite{gB00} it was proved that there exists an
$FB$-domain in $\C^n$ which is contained in the union of
$\Delta(r^2)\times\Delta^{n-1}(r)$ and the set
$S_1=\{z=(z_1,\dots, z_n):\ |z_1|\geq r^2-3r+||(z_2,\dots,
z_n)||_{\infty}\}$, for some $r>4$. This $FB$-domain is a basin of
attraction at 0 of a polynomial map that fixes the origin.
Therefore 0 is in the $FB$-domain. By using rotations, we deduce
that there exists an $FB$-domain which contains the origin and is
contained in the union of
$\Delta^{k-1}(r)\times\Delta(r^2)\times\Delta^{n-k}(r)$ and the
set $S_k=\{z=(z_1,\dots, z_n):\ |z_k|\geq r^2-3r+||(z_1,\dots,
\hat{z}_k,\dots, z_n)||_{\infty}\}$ for some $r>4,\ \forall
k=1,\dots, n$, where $\hat{z}_k$ means that $z_k$ is omitted. It
follows that in $\C^n$ there are $n$ $FB$-domains whose
intersection is non-empty and bounded. By post-composing the
corresponding $FB$-maps with contractions, and using translations,
we conclude that intersections of $FB$-domains form a base of
neighborhoods at each point of $\C^n$.

Now, under the assumptions of Theorem~\ref{T:1}, we can prove that
$\psi$ is continuous. Let $f_1,\dots, f_n$
be $FB$-maps as above so that the intersection of their images is
bounded.  Then, by assumption, $g_i=\psi\circ f_i\circ\psi^{-1},\
i=1,\dots, n$ are holomorphic maps in $\mathcal M$. Moreover,
$$
\psi(f_1(\C^n)\cap\dots \cap f_n(\C^n))=\psi(f_1(\C^n))\cap\dots
\cap\psi(f(\C^n))
$$
and since each $g_i$ is an injective holomorphic map,
$\psi(f_i(\C^n))=g_i(\psi(\C^n))=g_i({\mathcal M})$ is an open set. It
follows that $\psi$ is an open map. Using this plus the fact that
$\psi$ is a bijection of $\C^n$ onto a manifold, we see that
$\psi^{-1}(K)$ is compact for each compact $K \subset {\mathcal M}$.
With a standard argument, we conclude that $\psi$ is a
homeomorphism. In particular~\cite{wH48}, the dimension of $\mathcal M$
must be equal to $n$.

Note that \cite{dV00} implies that a Stein manifold $\mathcal M$ with
the (volume) density property has an
injective holomorphic map $F:\ {\mathcal M}\to{\mathcal M}$ with
$F({\mathcal M})\neq{\mathcal M}$. Since our proof of the continuity
of $\psi$ in Theorem~\ref{T:1} is based on the existence of special
maps of this form in $\C^n$, it is an interesting open
question whether such $\mathcal M$ can be shown to have a base of
neighborhoods given by finite intersections of injective images of
$\mathcal M$. If so, then it should be possible to replace $\C^n$ in
Theorem~\ref{T:1} by any manifold with these properties.

\section{Local linearization of maps}\label{S:LL}

Having shown that $\psi$ is continuous, we proceed as in~\cite{sM02}
to prove that $\psi$ is biholomorphic or antibiholomorphic using a
simultaneous linearization of certain commuting maps.
Let $a\in \C^n$ be an arbitrary point, and $b=\psi(a)$.
It is enough to show that $\psi$
is biholomophic or antibiholomorphic in a neighborhood of $a$.

A set ${\mathcal P}=\{p_i\}_1^n$ will be called a {\emph{system of
projections}} at $o$ in a complex manifold ${\mathcal N},\ o\in {\mathcal N}$,
if it consists of holomorphic maps in $E({\mathcal N})$
that fix $o$, and satisfy:
\begin{enumerate}
\item $p_i\neq o,\ \  \forall i$;
\item $p_i^2=p_i,\ \  \forall i$;
\item $p_i\circ p_j=o,\ \  \forall i\neq j$,
\end{enumerate}
where $p_i^2=p_i\circ p_i$, $o$ in (1) and (3) stands for the
constant map sending $\mathcal N$ to $o\in\mathcal N$. Let $f$ be
a biholomorphic map of $\mathcal N$ onto itself, that commutes
with all maps of some system of projections ${\mathcal P}$ at $o$,
and fixes $o$. We also assume that for every neighborhood $U$ of
$o$, and every compact set $K$, there exists an iterate of $f$
that brings $K$ into $U$, i.e. there exists a positive integer $l$
such that $f^l(K)\subset U$. Such a map $f$ clearly exists if
${\mathcal N}=\C^n$, since we can take it to be a contraction at
$o$, and $\{p_i\}$ to be standard projections. Now we introduce a
subsemigroup $I_f$ of $E({\mathcal N})$, consisting of all maps
$h$ that satisfy all the properties that $f$ does, with the same
system of projections $\mathcal{P}$, and such that $h$ commutes
with $f$. For reasons that will be clear later, we call the triple
$\{f, {\mathcal P}, I_f \}$ a {\emph{linearizing triple}}. It is
immediate to verify that all properties listed for a linearizing
triple are preserved under conjugation by $\psi$, i.e. if  $\{f,
{\mathcal P}, I_f \}$ is a linearizing triple in $\C^n$ at $a$,
then $\{g, {\mathcal Q}, I_g\}$ is a linearizing triple in
$\mathcal M$ at $b$, where $g=\psi\circ f\circ\psi^{-1},\
{\mathcal Q}=\psi\circ{\mathcal P}\circ\psi^{-1}$.

We note here that in general it is impossible to
linearize a holomorphic map in a neighborhood of its attracting
fixed point due to the presence of resonances among the
eigenvalues of its linear part~\cite{vA88}, \cite{jR88}. However, as
seen in the following proposition, under the
assumption that $h\circ p_i=p_i\circ h,\ \forall i=1,\dots, n$, the local
linearization of $h\in I_f$ is possible.

\begin{proposition}\label{P:Lin}
For every linearizing triple $\{f, {\mathcal P}, I_f\}$ in a complex manifold
$\mathcal N$ at $o$, there exists a biholomorphic map $\theta$
from a neighborhood of $o$ onto a neighborhood of the origin in $\C^n$,
such that for every $h\in I_f$, in some neighborhood of $o$,
\begin{equation}\label{E:L}
\theta\circ h=\Lambda_h\circ\theta,
\end{equation}
\begin{equation}\label{E:P}
\theta\circ p_i=P_i\circ\theta,\ \ \forall i,
\end{equation}
where $\Lambda_h$ is a diagonal linear map $(z_1,\dots,
z_n)\mapsto(\lambda_1 z_1,\dots, \lambda_n z_n)$, $\lambda_i, \
i=1,\dots, n$ satisfy $0<|\lambda_i|<1$, and are eigenvalues of
the linear part of $h$ at $o$, and $P_i$ is a diagonal matrix
similar to the linear part of $p_i$ at $o$.
\end{proposition}

The proof of this proposition follows the same arguments
as in~\cite{sM02}, and therefore we give only an outline here.
Because of the property that for every arbitrary compact set and every
neighborhood of
$o$, some iterate of $f$ brings the compact set into that
neighborhood, it follows that the eigenvalues of the linear part
of $f$ at $o$ are smaller than 1 in absolute value. Using the fact that
projections are locally linearizable~\cite{sK98}, and the commutativity
relations $h\circ p_i=p_i\circ h,\ \forall i$, the problem about local
linearization reduces to the one-dimensional Schr\"{o}der equation, which
is solved~\cite{aE93}. That all maps $h$ are linearized by the same
biholomorphic map $\theta$ follows from the
uniqueness of the solution to the Schr\"{o}der equation, and
the commutativity relations between $f,\ h$, and $p_i$.

We see in the following lemma that all invertible diagonal linear
maps whose entries are smaller than 1 in absolute value appear
in~(\ref{E:L}).

\begin{lemma}\label{L:B}
The map $\theta$ extends to a biholomorphic map of $\mathcal N$
onto $\C^n$. Moreover, if $\Lambda$ is a diagonal linear map
\begin{equation}
(z_1,\dots, z_n)\mapsto(\lambda_1 z_1,\dots, \lambda_n z_n),\ \
0<|\lambda_i|<1,\ i=1,\dots, n,\notag
\end{equation}
then there exists $h\in I_f$, such that
\begin{equation}\label{E:B}
\theta\circ h=\Lambda\circ\theta.
\end{equation}
\end{lemma}
\emph{Proof.} First we show that the map $\theta$ extends to a
biholomorphic map on the whole $\mathcal N$. This can be seen by using the
formula
\begin{equation}\label{E:T}
\theta=\Lambda_f^{-k}\circ\theta\circ f^k,\ \ k=1, 2, \dots
\end{equation}
Because of the property that for every compact subset $K$ of
$\mathcal N $ and every neighborhood $U$ of $o$, some iterate of
$f$ brings $K$ into $U$, it follows from~(\ref{E:T}) that $\theta$
can be extended to larger and larger sets, until its domain fills
the whole $\mathcal N$. Since $f$ and $\Lambda_f$ are
biholomorphisms, $\theta$ is injective, and hence a biholomorphism
on $\mathcal N$. The inverse of $\theta$ has a representation
similar to~(\ref{E:T}), and therefore $\theta$ is onto.

Consider a map $h=\theta^{-1}\circ\Lambda\circ\theta\in E(\mathcal
N)$. It is a biholomorphism of $\mathcal N$ onto itself, and it
commutes with every $p_i$, which follows from~(\ref{E:P}). Since
all entries of $\Lambda$ are less than 1 in absolute value, it is
clear that for every compact set $K$ and a neighborhood $U$ of
$o$, some iterate of $h$ brings $K$ into $U$. Using~(\ref{E:L}),
we conclude that $h$ commutes with $f$, and thus it belongs to
$I_f$. $\Box$

\section{Matrix equation}\label{S:ME}

Using the results of the previous section, we convert the statement of
Theorem~\ref{T:1} to a linearized version, thus reducing the problem
to determining the exact form of the solution of a matrix
equation ((\ref{E:RM}) below).
By finding this solution, we obtain an explicit expression for a map
$L$ defined below, which is conjugate to $\psi$ via biholomorphic
maps. This, with some more effort, will lead us to the proof that
$\psi$ is either biholomorphic, or antibiholomorphic.

We denote by ${\mathcal D}_0$ the set of invertible diagonal
$n\times n$ matrices whose entries are less than 1 in absolute
value, and we denote by ${\mathcal D}_n$
the set of all diagonal $n\times n$ matrices.
We identify ${\mathcal D}_0$ with the set of
diagonal linear maps, and ${\mathcal D}_n$ with a multiplicative
semigroup $\C^n$ in the obvious way, and consider a topology on
${\mathcal D}_n$ induced by the standard topology on $\C^n$.

In the previous section, we showed that if $\{f, {\mathcal P},
I_f\}$ is a linearizing triple in $\C^n$ at $a$, then $\theta:
\C^n \rightarrow \C^n$ conjugates $I_f$ to the set of diagonal
linear maps, which is isomorphic to ${\mathcal D}_0$. Similarly,
for a linearizing triple $\{g, {\mathcal Q}, I_g\}$ at
$b\in\mathcal M$, where $g=\psi\circ f\circ\psi^{-1},\ {\mathcal
Q}=\psi\circ{\mathcal P}\circ\psi^{-1}$, $I_g$ is conjugated by a
biholomorphic map $\eta: {\mathcal M} \rightarrow \C^n$ to
${\mathcal D}_0$.

We define a homeomorphism $L$ on $\C^n$ by
$L=\eta\circ\psi\circ\theta^{-1}$. For every $\Lambda$ in
${\mathcal D}_0$ we have
\begin{align}
L\circ\Lambda\circ L^{-1}&=\eta\circ\psi\circ\theta^{-1}
\circ\Lambda\circ\theta\circ\psi^{-1}\circ\eta^{-1}\notag \\
&=\eta\circ\psi\circ h\circ\psi^{-1}\circ\eta^{-1}=\eta\circ
j\circ\eta^{-1} =M\notag,
\end{align}
where $h=\theta^{-1}\circ\Lambda\circ\theta\in I_f;\ j=\psi\circ
h\circ\psi^{-1},\ M=\eta\circ j\circ\eta^{-1},\ M\in {\mathcal
D}_0$. Therefore the conjugation by $L$ defines an injective map $R$ from
${\mathcal D}_0$ to ${\mathcal D}_0$,
$R(\Lambda)=L\circ\Lambda\circ L^{-1}$, which is trivially
multiplicative, i.e.
$R(\Lambda'\Lambda'')=R(\Lambda')R(\Lambda''),\ \ \Lambda',\
\Lambda'', \Lambda'\Lambda''\in {\mathcal D}_0.$ Since $R$ is
continuous, it extends to a multiplicative map, which will also be
denoted by $R$ for convenience, from the subset
$\overline{\mathcal D}_0$ of ${\mathcal D}_n$ that consists of all
matrices in ${\mathcal D}_n$ with entries less than or equal to 1
in absolute value, into itself. Indeed, for every matrix $\Gamma$
in $\overline{\mathcal D}_0$, the image $R(\Gamma)$ also belongs
to $\overline{\mathcal D}_0$, which follows from the continuation
process. Now we extend $R$ to all of ${\mathcal D}_n$ as follows.
Let $\Gamma$ be an arbitrary matrix in ${\mathcal D}_n$. We choose
$\Lambda={\rm{diag}}(\lambda_1,\dots, \lambda_n)$ in ${\mathcal
D}_0$ such that $\sum_{i=1}^n|\lambda_i|\leq1$ and
$\Gamma\Lambda\in\overline{\mathcal D}_0$. Define
\begin{equation}\label{E:DR}
R(\Gamma)=R(\Gamma\Lambda)R(\Lambda)^{-1}\notag.
\end{equation}
The extended map $R$ is well defined. Indeed, if $\Lambda'$ is a
different matrix with the same properties as $\Lambda$, then
$R(\Gamma\Lambda)R(\Lambda')=R(\Gamma\Lambda')R(\Lambda)$, and the
conclusion follows from the commutativity relations for diagonal matrices.
The map $R$ is clearly
injective, and
multiplicative,
\begin{equation}\label{E:RM}
R(\Lambda'\Lambda'')=R(\Lambda')R(\Lambda''),\ \ \Lambda',\
\Lambda''\in {\mathcal D}_n.
\end{equation}

We denote by $\delta_i$ the diagonal $n\times n$ matrix which has
1 as its $ii$'th entry and all other entries 0. The system
$\{\delta_i\}_{i=1}^n$ is clearly the only one in ${\mathcal D}_n$
which satisfies $\delta_i\neq0,\ \delta_i^2=\delta_i,\
\delta_i\delta_j=0,\ \forall i\neq j$. Therefore, injectivity of
$R$ and~(\ref{E:RM}) imply that $R(\delta_i)=\delta_j,\ \forall
i$, where $j=j(i)$ is a permutation. In particular,
\begin{equation}\label{E:DE}
R(\delta_i\Lambda)=\delta_jR(\Lambda).
\end{equation}
If we denote the $jj$'th entry of the diagonal matrix $R(\Lambda)$
by $r_j(\lambda_1,\dots, \lambda_n)$, then~(\ref{E:DE}) implies
that $r_j$ depends on $\lambda_i$ only. For convenience, we write
$r_j(\lambda_1,\dots, \lambda_n)=r_j(\lambda_i)$. We can
rewrite~(\ref{E:RM}) as
\begin{equation}\label{E:Er}
r_j(\lambda_i'\lambda_i'')=r_j(\lambda_i')r_j(\lambda_i''),\ \
i=1,\dots, n,\ j=j(i).
\end{equation}
As in~\cite{aE93}, for every $j=j(i)$, the equation~(\ref{E:Er}) has
either the constant solution $r_j(\lambda_i)=1$, or
\begin{equation}\label{E:EE}
r_j(\lambda_i)=\lambda_i^{\alpha_{ij}}\overline{\lambda}_i^{\beta_{ij}},\
\ \alpha_{ij}, \beta_{ij}\in\C,\ \alpha_{ij}-\beta_{ij}=\pm1,
\end{equation}
where the last relation between $\alpha_{ij}$ and $\beta_{ij}$ is forced by
the injectivity of the map $R$.
Using~(\ref{E:EE}), we can obtain an explicit expression for $L:$
\begin{align}
L(z_1,\dots, z_n)&= {\rm{diag}}(z_{i(1)}^{\alpha_{i(1)1}}
\overline{z}_{i(1)}^{\beta_{i(1)1}},\dots,
z_{i(n)}^{\alpha_{i(n)n}}
\overline{z}_{i(n)}^{\beta_{i(n)n}})L(1,\dots, 1)\notag\\
&=B(z_1^{\alpha_1}\overline{z}_1^{\beta_1},\dots,
z_n^{\alpha_n}\overline{z}_n^{\beta_n}),\label{E:LE}\ \
\alpha_i-\beta_i=\pm1,\ i=1,\dots, n,
\end{align}
where $i=i(j)$ is an inverse permutation to $j=j(i)$, and $B$ is a
constant matrix.

By definition, $\psi=\eta^{-1}\circ L\circ\theta$. From the
expression~(\ref{E:LE}) for $L$ we can conclude that $\psi$ is
$\R$-differentiable and non-degenerate in $\C^n\setminus
\theta^{-1}(A)$, where $A=\cup_{k=1}^n\{(z_1,\dots, z_n):\
z_k=0\}$. Since the set $\theta^{-1}(A)$ is analytic, and using
the standard continuation argument for holomorphic maps, we can
assume that the map $\psi$ is $\R$-differentiable and
non-degenerate everywhere in $\C^n$. However, this is possible if
and only if $\alpha_i+\beta_i=1,\ i=1,\dots, n$. Combining this
with the equation $\alpha_i-\beta_i=\pm1$, we deduce that either
$\alpha_i=1,\ \beta_i=0$, or $\alpha_i=0,\ \beta_i=1$.

It remains to show that either $\alpha_i=1,\ \forall i$, or
$\alpha_i=0,\ \forall i$. To get a contradiction, suppose that
\begin{equation}
L(z_1,\dots, z_n)=B(\dots, z_i,\dots, \overline{z}_j,\dots).\notag
\end{equation}
Then
\begin{equation}
L^{-1}(w_1,\dots, w_n)=(\dots, l_i(w_1,\dots, w_n),\dots,
l_j(\overline{w}_1,\dots, \overline{w}_n),\dots),\notag
\end{equation}
where $l_i,\ l_j$ are nonconstant, linear holomorphic functions. Let
$\theta=(\theta_1,\dots, \theta_n)$. We consider a map $h\in
E$ in the form
\begin{equation}
h=\theta^{-1}(\dots, \theta_i\theta_j,\dots,
\theta_j,\dots)\theta,\notag
\end{equation}
where $\theta_i\theta_j$ is the $i$'th coordinate, and $\theta_j$
is the $j$'th coordinate. Using the definition of $L$, we obtain
\begin{align}
\eta\circ\psi\circ h&\circ\psi^{-1}\circ\eta^{-1}=
L\circ\theta\circ h\circ\theta^{-1}\circ L^{-1}\notag\\
&=B'(\dots, l_i(w_1,\dots, w_n)l_j(\overline{w}_1,\dots,
\overline{w}_n),\dots, \overline{l_j(\overline{w}_1,\dots,
\overline{w}_n)},\dots)\notag
\end{align}
for some constant matrix $B'$. This map, and hence $\psi\circ
h\circ\psi^{-1}$ is not holomorphic though, which is a contradiction.
$\Box$

\section{Epimorphism between semigroups}\label{S:E}

In this section we give a proof of Theorem~\ref{T:2}.

For a complex manifold $\mathcal N$ we denote by $C({\mathcal N})$
the subsemigroup of $E({\mathcal N})$ consisting of constant maps.
If ${\mathcal N}=\C^n$, we denote $C=C(\C^n)$. In other words,
\begin{equation}\label{E:CM}
c\in C({\mathcal N}) {\text{ if and only if }} \forall f\in E({\mathcal N}),
\  c\circ
f=c.
\end{equation}
There is a natural one-to-one correspondence between the constant
maps in $E({\mathcal N})$ and points of $\mathcal N$:
for each $z\in {\mathcal N}$ there exists
$c_z$ that maps $\mathcal N$ to $z$, and conversely, for each $c\in
C({\mathcal N})$
there exists $z\in {\mathcal N}$, such that $c=c_z$.
\begin{lemma}\label{L:Bl}
Let ${\mathcal N}_1$ and ${\mathcal N}_2$ be complex manifolds,
with ${\mathcal N}_1$ being weakly doubly-transitive. Let $\Phi:\
E({\mathcal N}_1) \to E({\mathcal N}_2)$ be an epimorphism of
semigroups. Then there exists a bijective map $\Psi:\ {\mathcal
N}_1 \to {\mathcal N}_2$ such that
\begin{equation}\label{E:BEq}
\Phi(f)=\Psi\circ f\circ\Psi^{-1},\ \ \forall f\in E({\mathcal
N}_1).
\end{equation}
\end{lemma}
{\emph{Proof.}} Because of~(\ref{E:CM}), and the assumption that
$\Phi$ is an epimorphism, for every $c\in C({\mathcal N}_1)$ we
have that $\Phi(c)\in C({\mathcal N}_2)$. Now we can define a map
$\Psi:\ {\mathcal N}_1 \to {\mathcal N}_2$ as follows
\begin{equation}
\Psi(z)=w {\text{ if and only if }} \Phi(c_z)=c_w.\notag
\end{equation}
Let $f$ be arbitrary map in $E({\mathcal N}_1)$. Then
\begin{equation}\label{E:Ef}
f\circ c_z=c_{f(z)}.
\end{equation}
Applying $\Phi$ to both sides of~(\ref{E:Ef}), we obtain
\begin{equation}\label{E:Ep}
\Phi(f)\circ c_{\Psi(z)}=c_{\Psi(f(z))}, \notag
\end{equation}
which is equivalent to
\begin{equation}\label{E:Be}
\Phi(f)\circ\Psi=\Psi\circ f.
\end{equation}

Equation~(\ref{E:Be}) implies surjectivity of $\Psi$. Indeed,
since $\Phi$ is an epimorphism, for every $w\in {\mathcal N}_2$,
there exists $f\in E({\mathcal N}_1)$, such that $\Phi(f)=c_w$.
Therefore, by~(\ref{E:Be}), $\Psi\circ f(z)= w,\ \forall
z\in{\mathcal N}_1$, which implies that $\Psi$ is onto.

We prove that $\Psi$ is injective by showing that for every $w\in
{\mathcal N}_2$ the full preimage $S_w=\Psi^{-1}(w)$ consists of
one point. Assume, by contradiction, that $S_w$ consists of more
than one point for some $w$. It cannot be all of ${\mathcal N}_1$,
since $\Psi$ is onto. Let $z_1$ be a point in $S_w$, such that in
arbitrary neighborhood of it there exist a point in ${\mathcal
N}_1\setminus S_w$. Let $z_2$ be arbitrary point in $S_w$,
different from $z_1$. From our assumption that ${\mathcal N}_1$ is
weakly doubly-transitive, it follows that there exists $h\in
E({\mathcal N}_1)$, such that $h(z_1)=z_1\in S_w$, and
$h(z_2)\notin S_w$. Evaluating $\Phi(h)$ at $w$, and
applying~(\ref{E:Be}) we have
\begin{align}
&\Phi(h)(w)=\Phi(h)\circ\Psi(z_1)=\Psi\circ h(z_1)=\Psi(z_1)=w,\notag\\
&\Phi(h)(w)=\Phi(h)\circ\Psi(z_2)=\Psi\circ h(z_2)\neq w,\notag
\end{align}
which is a contradiction. Thus we proved that $\Psi$ is a
bijection, and the equation~(\ref{E:BEq}) follows
from~(\ref{E:Be}). $\Box$

The first part of Theorem~\ref{T:2} now follows from
Lemma~\ref{L:Bl} and Theorem~\ref{T:1}, if we choose ${\mathcal
N}_1=\C^n$, and ${\mathcal N}_2={\mathcal M}$. The second part
follows if we take ${\mathcal N}_1={\mathcal M},\ {\mathcal
N}_2=\C^n$, and observe that equation~(\ref{E:BEq}) implies that
$\Phi$ is an isomorphism. $\Box$

\section{Isomorphic semigroups for inequivalent manifolds}\label{S:I}

In this section we prove Theorem~\ref{T:3}. We construct the
domains $D_1$ and $D_2$ by taking direct sums of $n$ copies of
domains as in Hinkkanen~\cite{aH92}. From \cite{aH92}, we know
that there exist unbounded domains $U_1, U_2$ in $\C$ such that
$U_1$ is neither conformally nor anticonformally equivalent to
$U_2$, and such that $E(U_1)$, and $E(U_2)$ are isomorphic and
consist of the constants plus the identity. One such choice of
domains is given by $U_1=\C\setminus\{0, 1, 2\}$, and
$U_2=\C\setminus\{0, 1, 2, \dots\}$. We set $D_1 = U_1 \times
\cdots \times U_1$, $D_2 = U_2 \times \cdots \times U_2$, and
verify that for these domains the conclusion of Theorem~\ref{T:3}
holds.

Let $F\in E(D_m),\ m=1,2$. Then each component $f_j$ of $F$ maps
$D_m$ holomorphically into $U_m$. Therefore, by the choice of
$U_m$, if we fix all $z_k,\ k=1,\dots, n, k\neq i$, then the
induced map $g_j(z_i)$ is in $E(U_m)$, hence is either a constant
map or the identity. Since $f_j$ is a continuous function in a
domain, which is a direct sum of domains in $\C$, we conclude that
it is identically equal to either a constant, or $z_i$ for some
$i=1,\dots, n$. Using this description of the elements in
$E(D_m)$, we can easily show that $E(D_1)$ and $E(D_2)$ are
isomorphic. Let $\xi$ be a bijective map from $U_1$ onto $U_2$. If
$F$ is an endomorphism of $D_1$, whose components are $f_1,\dots,
f_n$, then we set $\phi(F)$ to be an endomorphism of $D_2$, whose
$j$'th component is $z_i$ if $f_j=z_i$, and $\xi(c)$ if $f_j=c$, a
constant map. It is a simple matter to verify that the map $\phi$,
so defined, is an isomorphism of semigroups.

To show that $D_1$ and $D_2$ are not biholomorphically or
antibiholomorphically equivalent, we argue by contradiction.
Suppose first that there exists a biholomorphic map $F$ from $D_1$
onto $D_2$. Let $g$ be a non-constant restriction of a component
of $F$ to a coordinate axis. Such a component exists, since
otherwise the map $F$ would be constant. Since $g$ omits more than
two points, each of the points $0, 1, 2, \infty$ must be a
removable singularity or a pole. Therefore, $g$ extends to a
rational map. But this is a contradiction, because $g$ omits
infinitely many points. Similarly, we arrive at a contradiction by
assuming that there exists an antibiholomorphic map from $D_1$
onto $D_2$, and applying the same argument to a conjugate map.
$\Box$

\end{document}